\pdfoutput=1
\documentclass[english,11pt]{article}
\usepackage[T1]{fontenc}
\usepackage[latin9]{inputenc}
\usepackage{geometry}
\usepackage{graphicx}
\usepackage{amssymb}
\usepackage[usenames,dvipsnames]{color}
\usepackage{amsmath}

\RequirePackage[colorlinks,citecolor=blue,urlcolor=blue]{hyperref}

\makeatletter

\providecommand{\tabularnewline}{\\}

\@ifundefined{showcaptionsetup}{}{%
 \PassOptionsToPackage{caption=false}{subfig}}
\usepackage{subfig}
\makeatother

\usepackage{babel}

\usepackage[blocks]{authblk}

\usepackage{mdwlist}

\usepackage{enumitem}

\definecolor{olivegreen}{RGB}{34,139,34}

\newcommand{\rem}[1]{}

\newcommand{\R}{\mathbb{R}}

\newcommand{\E}{\mathbb{E}}

\newtheorem{proposition}{Proposition}
\newtheorem{theorem}{Theorem}

\def\Cplusplus{{\rm C\raise.3ex\hbox{\small++ }}}

\usepackage{array}
\usepackage{url}

\newcommand\Mark[1]{\textsuperscript#1}

\begin{document}

\begingroup
\centering
{\LARGE Technical Note for Discrete-Time Diffusion Approximations 
Motivated from Hospital Inpatient Flow Management \\[1.5em]
\large J. G. Dai\Mark{1} and Pengyi Shi\Mark{2}} \\[1em]
\begin{tabular}{*{1}{>{\centering}p{.85\textwidth}}}
\Mark{1}School of Operations Research and Information Engineering, Cornell University\tabularnewline
\url{jd694@cornell.edu}\tabularnewline
\Mark{2}Krannert School of Management, Purdue University  \tabularnewline
\url{shi178@purdue.edu}
\end{tabular}\par\par
\endgroup



\begin{abstract}
This note details the development of a discrete-time diffusion process 
to approximate the midnight customer count process 
in a $M_\textrm{per}/\textrm{Geo}_\textrm{2timeScale}/N$ system. 
We prove a limit theorem that supports this diffusion approximation, 
and discuss two methods to compute the stationary distribution of this discrete-time diffusion process.
\end{abstract}

\bigskip

The $M_\textrm{per}/\textrm{Geo}_\textrm{2timeScale}/N$ is single-pool queueing system
with a periodic Poisson arrival process and a \emph{two-time-scale} service time feature. 
This queueing system is motivated to study hospital inpatient flow management 
and is introduced in~\cite{Dai-Shi2012} (which we will refer to as the ``main paper''). 
To analyze this system, a critical step is to obtain the stationary distribution $\pi$
for the midnight count process $\{X_k ,k=0, 1, \dots\}$, 
where $X_k$ denotes the number of customers in the system at the midnight of day $k$, 
including both customers in service and those waiting in the buffer. 
To efficiently compute $\pi$, especially when the number of servers $N$ is large and the utilization is close to 1, 
we develop a discrete-time diffusion process $\{X^*_k, k=0, 1, \dots\}$ to
approximate the midnight count process, 
and use the stationary distribution of $\{X^*_k\}$ to approximate $\pi$. 

In this note, we first prove a limit theorem in Section~\ref{sec:diffusion-limit-single-cluster}
that supports this diffusion approximation. 
Then, in Section~\ref{sec:compute-stationary-for-diffusion} we discuss two methods
to numerically compute or approximate the stationary distribution of the discrete-time diffusion process $\{X^*_k\}$.  
Finally, in Section~\ref{sec:numerical-results-for-diffusion}, 
we show the accuracy of these two methods in approximating $\pi$ via numerical experiments.

\section{Diffusion limits for the single-pool model}
\label{sec:diffusion-limit-single-cluster}

Section 4.3 of~\cite{Dai-Shi2012} proposes 
a discrete-time diffusion process to approximate the midnight count process. 
This approximation is motivated by a limit theorem that shows the convergence of stochastic processes.
In this section, we prove this limit theorem. 

Instead of fixing the number of servers $N$,
we consider a sequence of 
$M_{\text{peri}}/\text{Geo}_{\text{2timeScale}}/N$ systems indexed by $N$,
i.e., a sequence of the single-pool models described in the main paper~\cite{Dai-Shi2012}. 
Let $\Lambda^N$ be the daily arrival rate of the $N$th system. 
Let $m = 1/\mu$, the mean LOS, be fixed and $\rho^N=(\Lambda^N m)/N$ be the
traffic intensity of the $N$th system. 
We assume that
\begin{equation}
\lim_{N\rightarrow \infty} \Lambda^N/N =\Lambda^*,
\textrm{ and }
\lim_{N\rightarrow \infty} \sqrt{N}(1-\rho^N)=\beta^* \textrm{ for some } \beta^*>0.
\label{eq:heavy-traffic-condition}
\end{equation}
Analogous to the conventional many-server queues that model customer
call centers~\cite{gkm03}, we call Condition~(\ref{eq:heavy-traffic-condition}) 
the \emph{Quality- and Efficiency-Driven} (QED) condition.  

We use $X^N_k$ to denote the midnight customer count at the midnight of day $k$ in the $N$th system.
We consider the \emph{diffusion-scaled} midnight customer count processes
$\tilde{X}^N = \{\tilde{X}_k^N: k=0, 1, 2, \ldots\}$ for the sequence of singe-pool systems, 
where for a given $k$, $\tilde{X}_k^N$ is defined as 
\begin{equation}
\tilde{X}_k^N = \frac{X_k^N-N}{\sqrt{N}} . 
\label{eq:diffusion-scale-process-define}
\end{equation}
Adapting the derivations in the main paper, 
we can show that $\tilde{X}_k^N$ satisfies the following relationship:  
\begin{equation}
     \tilde{X}^N_k = \tilde{Y}^N_k + \mu \sum_{i=0}^{k-1}(\tilde{X}^N_i)^- , 
     \quad k=0, 1, 2 \ldots, 
\label{eq:pre-limit-process-linear-form}     
\end{equation}
where  
\[
\tilde{Y}^N_k = \tilde{X}^N_0+\frac{1}{\sqrt{N}} \left(A^N_{(0,k]} - k\Lambda^N \right)
-\frac{1}{\sqrt{N}} \left( D^N_{(0,k]} - \mu (Z^N_0+\ldots+ Z^N_{k-1}) \right)
+k\sqrt{N}\mu(\rho^N-1),
\]
$A^N_{(0,k]} = \sum_{i=0}^{k-1} A^N_i $ and $D^N_{(0,k]} = \sum_{i=0}^{k-1} D^N_i$ 
are the cumulative number of arrivals and departures from 0 until the midnight (zero hour) of day $k$ 
in the $N$th system, respectively, 
and $Z^N_i = \min (X^N_i, N)$ is the number of busy servers at the midnight of day $i$.
We assume the initial condition
\begin{equation}
\tilde{X}^N_0 \Rightarrow {X}^*_0 \textrm{ as } N\rightarrow \infty,
\label{eq:initial-condition}
\end{equation}
where $\Rightarrow$ denotes convergence in distribution.
Under the many-server heavy-traffic framework (e.g., see~\cite{DaiHeTezcan10}),
we prove the following limit theorem: 
\begin{theorem}\label{conj:htlSingle}
Consider a sequence of $M_{\text{peri}}/\text{Geo}_{\text{2timeScale}}/N$
single-pool systems that satisfies (\ref{eq:heavy-traffic-condition}) and (\ref{eq:initial-condition}).
For any positive integer $K \in \mathbb{Z}_+$, 
$\tilde{X}^N \Rightarrow X^{\ddagger} $ on the compact set $\left[ 0,K \right]$
as $N\rightarrow \infty$, i.e., 
\begin{equation}
\left( \tilde{X}^N_0, \tilde{X}^N_1, \dots, \tilde{X}^N_K \right) \Rightarrow
\left( X^{\ddagger}_0, X^{\ddagger}_1, \dots, X^{\ddagger}_K \right) \textrm{ as } N \rightarrow \infty.
\label{eq:convergence-driven-process-compact-set}
\end{equation}   
The discrete-time limit process $X^{\ddagger} = \{ X^{\ddagger}_k, k=0, 1, \dots \}$ satisfies
\begin{equation}
    X^{\ddagger}_k = Y^{\ddagger}_k + \mu \sum_{i=0}^{k-1}(X^{\ddagger}_i)^-, \quad k=0, 1, \dots, 
\label{eq:diffusion-limit-theorem}
\end{equation}
where $Y^{\ddagger}_k=Y^{\ddagger} (k)$ for $k=0, 1, \ldots,$ is an embedding of the Brownian motion 
$\{Y^{\ddagger} (t), t\ge 0\}$ which starts from $X_0^{\ddagger}$ 
and has mean $-\mu\beta$ and variance $\Lambda^*+\mu(1-\mu)$. 
\end{theorem}
In this limit theorem, we deliberately use the superscript $\ddagger$ to differentiate 
the limit process $X^{\ddagger}$ (and the associated process $Y^{\ddagger}$) 
from the diffusion \emph{approximation}
$X^*$ (and the associated process $Y^{*}$) introduced in Section~4.3 of the main paper.

The key step of the proof for Theorem~\ref{conj:htlSingle} is to show that 
$\{\tilde{Y}^N_k, k=0,1,\dots\}$ converges to $\{Y^{\ddagger}_k, k=0,1,\dots\}$
on any given compact set $\left[ 0,K \right]$, 
or equivalently, 
\begin{equation}
\left( \tilde{Y}^N_0, \tilde{Y}^N_1, \dots, \tilde{Y}^N_K \right) \Rightarrow
\left( Y^{\ddagger}_0, Y^{\ddagger}_1, \dots, Y^{\ddagger}_K \right) \textrm{ as } N \rightarrow \infty.
\label{eq:convergence-driven-process-compact-set}
\end{equation}
Then, the convergence of $\tilde{X}^N$ to $X^{\ddagger}$ naturally follows
because of the linear forms in~(\ref{eq:pre-limit-process-linear-form}) and~(\ref{eq:diffusion-limit-theorem}). 
To prove~(\ref{eq:convergence-driven-process-compact-set}),
we first prove the convergence of the diffusion-scaled arrival processes in Section~\ref{sec:converge-arr-process}, 
and then the convergence of the discharge processes in Section~\ref{sec:converge-dis-process}.


\subsection{Arrival process}
\label{sec:converge-arr-process}

For the $N$th system, let $\tilde{E}^N_{k} = \frac{1}{\sqrt{N}} \left(A^N_{(0,k]} - k\Lambda^N \right)$.
We also introduce a continuous-time process $\{\tilde{E}^N (t), \ t\geq 0 \}$ defined as  
\begin{equation}
\tilde{E}^N (t) =  \frac{1}{\sqrt{N}} \left( E^N (t) - \Lambda^N t \right), 
\end{equation}
where $E^N (\cdot)$ represents a Poisson process with rate $\Lambda^N$. 
It is easy to verify that $\{ \tilde{E}^N_{k} \}$ is an embedding of $\tilde{E}^N (\cdot)$, i.e., 
\[
\tilde{E}^N_{k} = \tilde{E}^N (k), \ k=0, 1, \dots.
\]
Following a standard functional central limit theorem argument, we can show that
\begin{equation}
\tilde{E}^N (\cdot) \Rightarrow E^{\ddagger} (\cdot) 
\label{eq:E(t)-spaceD}
\end{equation}
in space $\mathbb{D}$ endowed with the Skorohod $J_1$ topology, 
where $E^{\ddagger} (\cdot)$ is a Brownian motion with drift 0 and variance $\Lambda^*$.
Because the convergence of stochastic processes implies 
the convergence of any finite-dimensional joint distributions, we then naturally have
\begin{equation}
\left( \tilde{E}^N_0, \tilde{E}^N_1, \dots, \tilde{E}^N_K \right) \Rightarrow
\left( E^{\ddagger}_0, E^{\ddagger}_1, \dots, E^{\ddagger}_K \right) \textrm{ as } N \rightarrow \infty, 
\label{eq:convergence-arrival-process-compact-set}
\end{equation}
where $E^{\ddagger}_k=E^{\ddagger} (k)$ is also an embedding of $E^{\ddagger} (\cdot)$.

\subsection{Discharge process}
\label{sec:converge-dis-process}

Now we consider the diffusion-scaled discharge processes. 
For the $N$th system, we introduce two discrete-time processes: 
\[
\check{D}^N_{k} = \frac{1}{\sqrt{N}} \left( D^N_{(0,k]} - \mu (Z^N_0+\ldots+ Z^N_{k-1}) \right), \ k=0, 1, \dots, 
\] 
and 
\[
\tilde{D}^N_{k} = \frac{1}{\sqrt{N}}\sum_{i=1}^{Z^N_0+\ldots+ Z^N_{k-1}} \left( \xi_i-\mu \right), \ k=0, 1, \dots, 
\]
where $\{\xi_{i}\}$ is a sequence of iid Bernoulli random variables with success probability $\mu$.
Recall that in Appendix C of the main paper, we establish a revised system 
which tosses coins for every customer in service at the midnight to determine the departures each day, 
and we have proved this revised system is equivalent to the original system in distribution. 
Using the revised system, we can show the above two discrete-time processes are equal in distribution, i.e., 
\[
(\check{D}^N_{0}, \check{D}^N_{1}, \dots ) =^d (\tilde{D}^N_{0}, \tilde{D}^N_{1}, \dots ). 
\]
Thus, it is sufficient to prove for any given $K \in \mathbb{Z}_+$, 
\begin{equation}
\left( \tilde{D}^N_0, \tilde{D}^N_1, \dots, \tilde{D}^N_K \right) \Rightarrow
\left( S^*_0, S^*_1, \dots, S^*_K \right) \textrm{ as } N \rightarrow \infty . 
\label{eq:convergence-discharge-process-compact-set}
\end{equation}
Here, $S^*_k=S^* (k)$ is an embedding of the Brownian motion $S^* (\cdot)$ 
with drift $0$ and variance $\mu(1-\mu)$.  

Let $\eta_{i} = \xi_i-\mu$, and $\{\eta_{i}\}$ forms a sequence of iid random variables 
with mean 0 and variance $\mu(1-\mu)$. 
We also define
\[
T^N_k = \sum_{j=0}^{k-1}{Z_j}, \quad \bar{T}^N_k = \frac{T^N_k}{N},
\]
and
\[
S_n = \sum_{i=1}^{n} \eta_{i}.
\]
Then, we can further rewrite $\tilde{D}^N_{k}$ as
\begin{eqnarray*}
\tilde{D}^N_{k} &=& \frac{1}{\sqrt{N}} \sum_{i=1}^{T^N_k}\eta_{i} = \frac{1}{\sqrt{N}} S_{\bar{T}^N_k N}. 
\end{eqnarray*}
Correspondingly, proving (\ref{eq:convergence-discharge-process-compact-set})
is equivalent to showing 
\begin{equation}
\left( \frac{1}{\sqrt{N}}S_{\bar{T}^N_0 N}, \frac{1}{\sqrt{N}}S_{\bar{T}^N_1 N}, \dots, \frac{1}{\sqrt{N}}S_{\bar{T}^N_K N} \right)
\Rightarrow
\left( S^*_0, S^*_1, \dots, S^*_K \right) \textrm{ as } N \rightarrow \infty.
\label{eq:convergence-discharge-process-compact-set-revised}
\end{equation}


To prove (\ref{eq:convergence-discharge-process-compact-set-revised}),
we introduce a continuous-time process $\{\tilde{S}^{N} (t), t\geq 0 \}$, where
\[
\tilde{S}^{N} (t) = \frac{1}{\sqrt{N}}S_{\lfloor t N \rfloor} \circ \bar{T}^N_{\lfloor t \rfloor}, \quad t\geq 0.
\]
In other words, $\tilde{S}^{N} (\cdot)$ is a composition of two continuous processes,
$\frac{1}{\sqrt{N}}S_{\lfloor \cdot N \rfloor}$ and $\bar{T}^N_{\lfloor \cdot \rfloor}$.
It is easy to verify that $\{ \tilde{D}^N_{k} \}$ is an embedding of $\tilde{S}^{N} (\cdot)$ 
because when $t=k$, $\bar{T}^N_k N = T^N_k $ is always an integer and 
\[
\tilde{D}^N_{k} = \frac{1}{\sqrt{N}}S_{\bar{T}^N_k N} = \tilde{S}^{N} (k). 
\] 
If we can show  
\begin{equation}
 \frac{1}{\sqrt{N}} S_{\lfloor \cdot N \rfloor} \Rightarrow S^* (\cdot) 
\label{eq:convergence-S-process-compact-set}
\end{equation}
in space $\mathbb{D}$ endowed with the Skorohod $J_1$ topology as well as 
\begin{equation}
\bar{T}^N_{\lfloor \cdot \rfloor} \rightarrow \bar{T}_{\lfloor \cdot \rfloor} \textrm{ in probability}
\label{eq:convergence-bar-T-process-compact-set}
\end{equation}
with $\bar{T}_{\lfloor t \rfloor} = \lfloor t \rfloor$,
then applying the random time change theorem, we can prove~(\ref{eq:convergence-discharge-process-compact-set-revised}).
The convergence in (\ref{eq:convergence-S-process-compact-set}) follows from the Donsker's theorem, 
and we focus on proving~(\ref{eq:convergence-bar-T-process-compact-set}) below.
It is sufficient to show for each $0\leq k \leq K$,
$Z_k^N / N \rightarrow 1$ almost surely, which we prove with induction.

We first rewrite the system equation under the fluid scaling: 
\begin{eqnarray}
\bar{X}^N_{k} &=& \bar{Y}^N_k + \sum_{i=0}^{k-1}(\bar{X}_{i}^{N})^-, \end{eqnarray}
where
\[
\bar{X}^N_{k} = \frac{X^N_k-N}{N},
\]
and
\[
\bar{Y}^N_k=\bar{X}^N_0+\frac{\sum_{i=1}^{k} (A^N_{i-1}-\Lambda^N)}{N}
                 -\frac{\sum_{i=1}^{T^N_k}\eta_{i}}{N}
		         +\frac{k(\rho^N-1)}{\sqrt{N}}.
\]
Assume that $X^N_0=N$, then $\bar{X}^N_0=0$ and $Z^N_0=N$
(so $\bar{X}^N_0 \rightarrow 0$ is trivial). 
\begin{itemize}

\item When $k=1$, we have
\[
\bar{Y}^N_1=\bar{X}^N_0+\frac{A^N_0 - \Lambda^N }{N}
                 -\frac{\sum_{i=1}^{N}\eta_{i}}{N} 
		+\frac{(\rho^N-1)}{\sqrt{N}}.
\]
Recall that $(A^N_0 - \Lambda^N)$ and $\eta_i$ are centered random variables with mean 0.
By the Law of Large Numbers, it is obvious that
\[
\bar{Y}^N_1 \rightarrow 0 \quad a.s. \textrm{ when } N\rightarrow \infty.
\]
Thus, $\bar{X}^N_1 = \bar{Y}^N_1 \rightarrow 0 \  a.s.$,
and $Z^N_1/N \rightarrow 1 \ a.s.$.

\item Assume that at $k$, we have for all $0 \leq j \leq k$,
$\bar{X}^N_j  \rightarrow 0 \ a.s.$
and $Z^N_j/N \rightarrow 1 \ a.s.$.  Then for $k+1$,
we have $\bar{T}^N_{k+1} \rightarrow (k+1) \ a.s.$ and
\begin{eqnarray}
\bar{Y}^N_{k+1} & = & \bar{X}^N_0+\frac{\sum_{i=1}^{k+1} (A^N_{i-1}-\Lambda^N) }{N}
                 -\frac{\sum_{i=1}^{T^N_{k+1}}\eta_{i}}{N}
		+\frac{(k+1)(\rho^N-1)}{\sqrt{N}} \nonumber \\
		& = & \frac{\sum_{i=1}^{k+1} (A^N_{i-1}-\Lambda^N) }{N} 
                 -\frac{\sum_{i=1}^{T^N_{k+1}}\eta_{i}}{T^N_{k+1}} \cdot \frac{T^N_{k+1}}{N}
		+\frac{(k+1)(\rho^N-1)}{\sqrt{N}} \\
		& \rightarrow & 0 \quad a.s..
\end{eqnarray}

Then
\[
\bar{X}^N_{k+1} = \bar{Y}^N_{k+1} + \sum_{j=0}^{k}(\bar{X}_{j}^{N})^-
\rightarrow 0 \quad a.s., 
\]
which completes the proof of Theorem~\ref{conj:htlSingle}.

\end{itemize}


\section{Computing the stationary distribution of the discrete-time diffusion process}
\label{sec:compute-stationary-for-diffusion}

Motivated by the limit theorem proved in Section~\ref{sec:diffusion-limit-single-cluster}, 
Section~4.3 of the main paper~\cite{Dai-Shi2012} proposes a discrete-time diffusion process 
$\{X_k^*, k=0,1,\dots \}$ to approximate the original midnight count process $\{X_k, k=0, 1, \dots \}$. 
The dynamics of this approximation process follows:
\begin{equation}
    {X}^*_k = {Y}^*_k + \mu \sum_{i=0}^{k-1}({X}^*_i)^-, \quad k=0, 1, 2, \ldots,
\label{eq:diffusion-limit-midnight-count}
\end{equation}
where ${Y}^*_k={Y}^*(k)$ for $k=0, 1, 2, \ldots,$
and $\{{Y}^*(t), t\ge 0\}$ is a Brownian motion
with mean
\begin{equation}
\theta_N = \Lambda - N\mu = - N \mu (1- \rho)
\label{eq:mean-diffusion-midnight-count-process}
\end{equation}
and variance
\begin{equation}
\sigma^2_N = \Lambda + \rho N \mu(1-\mu) = \rho N \mu(2-\mu).
\label{eq:variance-diffusion-midnight-count-process}
\end{equation}
Note that (\ref{eq:mean-diffusion-midnight-count-process})
and (\ref{eq:variance-diffusion-midnight-count-process})
are different from the mean $-\mu\beta$ and variance $\Lambda^*+\mu(1-\mu)$
in Theorem~\ref{conj:htlSingle}, for two reasons: 
first, the process ${X}^*_{\cdot}$ and ${Y}^*_{\cdot}$ 
are diffusion \emph{approximations} instead of the limiting processes stated in Theorem~\ref{conj:htlSingle},
which is why the term $\rho$ appears in (\ref{eq:mean-diffusion-midnight-count-process})
and (\ref{eq:variance-diffusion-midnight-count-process});
second, the process ${X}^*_{\cdot}$ is to approximate the \emph{centered} midnight count process 
(defined as $\hat{X}_k = X_k - N$), not the diffusion-scaled version as in (\ref{eq:diffusion-scale-process-define}). 

In the next three subsections, we first specify the basic adjoint relationship (BAR) 
for this discrete-time diffusion process ${X}^*_{\cdot}$. 
Then, we discuss two ways to numerically calculate/approximate the stationary distribution of ${X}^*_{\cdot}$:
(i) a projection algorithm that numerically solves the BAR, and (ii) an approximate formula. 

\subsection{Basic adjoint relationship}

The state space of $\{X_k^*, k=0,1,\dots \}$ is  $\R$. 
One can check that $\{X_k^*, k=0,1,\dots \}$ is a discrete-time Markov process, since
\[
X^*_{k+1} - X^*_{k} = Y^*_{k+1} - Y^*_{k} + \mu{(X^*_k)}^-, \text{ for } k=0, 1, \ldots,
\]
and $\{Y^*_{k+1}-Y^*_k: k=0,1, \ldots \}$ is a sequence of iid normal r.v.\ with mean $\theta_N$ and variance $\sigma_N^2$.
The transition density of the Markov process is
\begin{equation}
p(x,y)=\mathbb{P}(X^*_{k+1}=y | X^*_k =x) =
\begin{cases}
 \phi_{\theta_N, \sigma^2_N} \bigl(y-x \bigr),  & \textrm{ when } x \geq 0, \\
\phi_{\theta_N,\sigma^2_N}\bigl(y-(1-\mu)x \bigr), & \textrm{ when } x < 0,
\end{cases}
\label{eq:transition-density-for-BAR}
\end{equation}
where $\phi_{\theta, \sigma^2}$
denotes the normal density function with mean $\theta$ and variance $\sigma^2$.
Let $C_b(\R)$ denote the set of bounded, continuous functions on $\R$. For each
$f\in C_b(\R)$, define
\[
\bold{P}f(x)=\int_\mathbb{R} p(x,y)f(y)dy \quad \text{ for each } x\in \R.
\]
One can check that $\mathbf{P}f\in C_b(\R)$. It follows that the stationary density $\pi(x)$ satisfies
\begin{equation}
\int_\mathbb{R} \bold{P}f(x) \pi(x) dx = \int_\mathbb{R} f(x) \pi(x) dx,  \quad \forall f \in C_b (\mathbb{R}),
\label{eq:BAR}
\end{equation}
or equivalently,
\begin{equation}
\int_\mathbb{R} \bold{L}f(x) \pi(x) dx = 0, \quad \forall f \in C_b (\mathbb{R}),
\label{eq:BAR-equivalent}
\end{equation}
with $\bold{L}f(x) = \bold{P}f(x) - f(x)$.  We call
(\ref{eq:BAR-equivalent}) the basic adjoint relationship (BAR) that
governs the stationary density of the discrete-time Markov process
$\{X_k^*, k=0,1,\dots \}$.


\subsection{A projection algorithm}
\label{sec:projection-algorithm}

The BAR (\ref{eq:BAR-equivalent}) is in the same format as (2.5) of~\cite{Dai-He2013}; 
the latter BAR is for the stationary density of a
(continuous-time) diffusion process. As such the algorithm developed in~\cite{Dai-He2013} 
can be applied to compute the stationary density $\pi^*$ 
of the discrete-time diffusion process $\{X_k^*, k=0,1,\dots \}$.
We outline the algorithm here, commenting on the differences when appropriate.

\subsubsection{Reference density and  the  space $L^2(\R, r)$}

To compute the stationary density $\pi^*$, we first need a reference density $r$
such that
\[
\int_\mathbb{R} r(x)  dx = 1.
\]
We use the approximate formula $\tilde{\pi}$ in Section~4.3.2 of the main paper 
(also see~\ref{eq:diffusion-density-revised} below) as the reference density $r$.

Next, we define the ratio function as:
\begin{equation}
q(x) = \frac{\pi^*(x)}{r(x)} \quad \text{ for } x\in \R.
\label{eq:ratio-function}
\end{equation}
With the given reference density $r$, if we can compute the ratio function $q$,
then we can compute the stationary density via
\[
\pi^*(x) = q(x)r(x) \textrm{ for } x \in \mathbb{R}.
\]

To compute $q$,  we plug (\ref{eq:ratio-function}) into (\ref{eq:BAR-equivalent}) and get
\begin{equation}
\int_\mathbb{R} \bold{L}f(x) q(x) r(x) dx = 0, \quad \forall f \in C_b (\mathbb{R}).
\label{eq:BAR-reference-density}
\end{equation}
Following the notation in~\cite{Dai-He2013}, we use $L^2(\mathbb{R}, r)$
to denote the space of all square-integrable functions on $\mathbb{R}$ with
respect to the measure that has density $r$.
Namely, $L^2(\mathbb{R}, r)$  is the set of measurable functions $f$ on $\R$ that satisfy
\begin{displaymath}
 \int_{\R} f^2(x) r(x) dx <\infty.
\end{displaymath}
We adopt the same inner product
on $L^2(\mathbb{R}, r)$ as in~\cite{Dai-He2013}, that is,
\begin{equation}
\langle f, \hat{f} \rangle = \int_\mathbb{R} f(x) \hat{f}(x) r(x) dx, \quad \textrm{ for } f, \hat{f} \in L^2(\mathbb{R}, r).
\label{eq:inner-product-L2-space}
\end{equation}
In (3.2) of~\cite{Dai-He2013}, the authors made an important assumption on the reference density.
Namely, they assumed that the reference density was chosen so that
\begin{equation}
  q\in L^2(\R, r).
\label{eq:L2}
\end{equation}
With our choice of the reference density $r$,  
we have been unable to verify that condition (\ref{eq:L2}) is satisfied. 
We leave it as a conjecture that condition (\ref{eq:L2}) is satisfied.
The remainder of this section assumes that the conjecture is true.

\subsubsection{Orthogonal projection}
\label{sec:othpro}

Note that the BAR (\ref{eq:BAR-reference-density}) is equivalent to
\[
\langle \bold{L}f, q \rangle = 0 \quad \text{ for each } f \in C_b (\mathbb{R}).
\]
Thus, $q$ satisfying the BAR is equivalent to $q$ being  \emph{orthogonal} to $\mathbf{L}f$ for each $f \in C_b (\mathbb{R})$.
We define a space $H$ as
\[
H = \textrm{ the closure of } \{ \bold{L}f: f \in C_b (\mathbb{R}) \},
\]
which is a subspace of $L^2(\mathbb{R}, r)$.
Therefore, $q$  satisfying the BAR  is equivalent to $q$ being orthogonal to space $H$.
Therefore, our task is to find a function $q$ that is orthogonal
to space $H$.  To do so, we consider a constant function $e$ with
$e(x) = 1$ for each $x \in \mathbb{R}$.
Since
\begin{equation}
  \label{eq:normalize}
\langle e, q \rangle = \int_\mathbb{R} e(x) q(x) r(x) dx = \int_\mathbb{R} \pi^*(x) dx = 1,
\end{equation}
one can check that $e \notin
H$ because otherwise $\langle e, q \rangle=0$, contradicting (\ref{eq:normalize}).
We use $\bar{e}$ to denote the projection of $e$ onto $H$.
Then, $e - \bar{e}\neq 0$ and it  must be orthogonal to $H$.
Once we have $\bar{e}$, we obtain the ratio function $q$ by
\[
q = \frac{e - \bar{e}}{|| e - \bar{e} ||^2},
\]
where $||\cdot ||$ is the induced norm from the inner product (\ref{eq:inner-product-L2-space})
with $|| f ||^2 = \langle f, f \rangle$ for $f \in L^2(\mathbb{R}, r)$.

\subsubsection{Finite-dimensional approximation}

The projection of $e$ onto $H$ can be expressed as
\begin{equation}
\bar{e} = \textrm{argmin}_{h\in H} || e - h ||.
\label{eq:projection-e-original}
\end{equation}
The space $H$ is linear and infinitely dimensional.
To compute the projection numerically, we use a finite-dimensional subspace $H_k$
to approximate $H$ and find the projection $\bar{e}_k$ of $e$ on $H_k$, namely,
\begin{equation}
\bar{e}_k = \textrm{argmin}_{h \in H_k} || e - h ||.
\label{eq:projection-e-finite}
\end{equation}

Let $C_k$ be a finite-dimensional, linear subspace of $C_b(\R)$. 
Then $H_k=\{\mathbf{L}f: f\in C_k\}$ is a finite-dimensional subspace of $H$.
Assume that  $\{f_i: i=1, 2, \dots, m \}\subset C_k$ is a basis of $C_k$.
Then, since the projection $\bar{e}_k \in H_k$,
it can be represented as a linear combination of  $\{\bold{L}f_i : i=1, 2, \dots, m \}$. That is,
\begin{equation}
\bar{e}_k = \sum_{i=1}^{m} \alpha_i \bold{L}f_i
\label{eq:bar-ek-linear-combination}
\end{equation}
where $\alpha_i \in \mathbb{R}$ for $i=1, 2, \dots, m $.

To compute the vector of  coefficients $\alpha=(\alpha_1, \ldots, \alpha_m)'$, we use the fact that $\langle e-\bar{e}_k, \bold{L}f_i  \rangle = 0$
for $i=1, 2, \dots, m$. Consequently, we obtain a system of linear equations
\begin{equation}
  \label{eq:Ab}
A \alpha = \beta,
\end{equation}
where $A_{ij} = \langle \bold{L}f_i , \bold{L}f_j  \rangle$
and $\beta_i = \langle e, \bold{L}f_i  \rangle$ for $i,j=1, \ldots, m$.
The matrix is symmetric, semi-positive definite, but can be singular.
Although the solution to the system of linear equations may not be unique,
projection $\bar{e}_k$ is unique. When $A$ is singular or nearly singular,
one can solve (\ref{eq:Ab}) by direct methods such as the QR decomposition and
the Cholesky decomposition or by iterative methods such as LSQR~\cite{paige1982lsqr}.
The Cholesky decomposition exploits the symmetric and semi-positive definite
properties of $A$ even when $A$ is singular~\cite{anderson1999lapack,Hammarling07lapack-stylecodes},
whereas QR decomposition does not. Unlike many other iterative methods,
LSQR can handle matrix $A$ when it is singular. LSQR does not exploit semi-positive definiteness.

Once we get the vector of coefficients
$\alpha=(\alpha_1, \ldots, \alpha_m)'$ by solving the system of linear equations (\ref{eq:Ab}),
we can compute $\bar{e}_k$ as in (\ref{eq:bar-ek-linear-combination}).
Eventually, we can approximately compute the stationary density $\pi^*$ as
\begin{equation}
\pi^* (x) \approx r(x) \frac{1- \bar{e}_k (x)}{|| e - \bar{e}_k ||^2} \quad \forall x \in \mathbb{R}.
\label{eq:final-approx-pi-with-bar-ek}
\end{equation}

\subsubsection{FEM implementation}
\label{sec:implem}

In our implementation, we use the finite element method (FEM) to
construct the approximate space $C_k$, following Section 3.3 of~\cite{Dai-He2013}.
The numerical results in this paper for approximating the stationary density $\pi$
with the projection algorithm all follow this FEM implementation.
In Proposition~3 of Dai and He~\cite{Dai-He2013}, they proved the convergence of using
(\ref{eq:final-approx-pi-with-bar-ek}) to approximate $\pi$ as $H_k\uparrow H$.
Their proof applies to our setting when (\ref{eq:L2}) is satisfied.


\subsection{Approximate formula for the stationary density}

In Section 4.3 of the main paper, the following formula $\tilde{\pi}$ is proposed as 
a proxy for the stationary density $\pi^*$ of the diffusion process $X^*$:
\begin{equation}
\tilde{\pi}(x)=
\begin{cases}
 \alpha_1 \exp \bigl( {2\theta_N x} / {\sigma^2_N} \bigr),  & x\geq 0; \\
 \alpha_2 \exp \bigl(- {(2\mu-\mu^2)(x -\theta_N/\mu)^2} / {2\sigma^2_N} \bigr), & x<0;
\end{cases}
\label{eq:diffusion-density-revised}
\end{equation}
where $\alpha_1$ and $\alpha_2$ are normalizing constants that make $\tilde{\pi}(x)$ continuous at zero
and $\int_\mathbb{R} \tilde{\pi} (x) dx = 1$. 

As mentioned in the main paper, 
the rationale of this approximate formula is 
based on the analogy between $\{{X}^*_k:k=0, 1, 2 \ldots\}$
and $\{\check{X}(t), t\ge 0\}$, where 
\begin{equation}
    \check{X}(t) = \check{Y}(t) + \mu \int_0^t (\check{X}(s))^- ds,
    \quad t\ge 0,
\label{eq:continous-diffusion-limit}
\end{equation}
and $\{\check{Y}(t), t\ge 0\}$ is a Brownian motion.
To get the stationary density of $\check{X}$,
Browne and Whitt~\cite{BrowneWhitt} have suggested that
since (i) $\check{X}$ is a Ornstein-Uhlenbeck (OU) process on $(-\infty, 0]$
and the stationary density of an OU process has a Gaussian form
and (ii) $\check{X}$ is a reflected Brownian motion (RBM) on $[0, \infty)$
and the stationary density of a RBM has an exponential form,
then the stationary density of $\check{X}$ can be obtained
by piecing together the Gaussian and exponential densities.

We use the same piecing technique in our setting.
Specifically, $\{{X}^*_k:k=0, 1, 2 \ldots, \}$ behaves as a discrete version of the OU process
on $(-\infty, 0]$ and as a reflected random walk on $[0, \infty)$.
We show in Proposition~\ref{pro:stationary_density_OU} below that 
the stationary density of the discrete-time OU (DOU) process also has a Gaussian form. 
For the reflected random walk,
existing research shows that it has an exponential tail~\cite{Kingman1963, Siegmund79, BlaGlyCDA2006}.
Therefore, we piece together a Gaussian density and an exponential density
and propose using (\ref{eq:diffusion-density-revised}) to approximate $\pi^*$.
In the next two subsections, 
we first prove that the stationary density of the discrete-time OU process has a Gaussian form, 
then we show the details of deriving formula (\ref{eq:diffusion-density-revised}).

\subsubsection{The stationary distribution of a discrete OU process}
\label{app:stationary-dist-for-dist-OU}

Similar to the continuous-time version of the Ornstein-Uhlenbeck process, 
we define its discrete-time version $\{{X}^{\textrm{DOU}}_k, k=0, 1, \dots\}$ as: 
\begin{equation}
X^{\textrm{DOU}}_k=Y^{\textrm{DOU}}_k-\mu\sum_{i=0}^{k-1} X^{\textrm{DOU}}_i, \quad k=0, 1, \dots 
\label{eq:discrete-OU}
\end{equation} 
where $\{Y^{\textrm{DOU}}_k :=  \sum_{i=0}^{k-1} \xi_i, k=0, 1, \dots \}$ is a Gaussian random walk, 
i.e., $\{\xi_i\}$ is a sequence of iid random variables following a normal distribution 
with mean $\theta$ and variance $\sigma^2$.

The following proposition says the stationary density for a discrete OU process
has the Gaussian form, which is consistent with that in a continuous-time OU process.
\begin{proposition} \label{pro:stationary_density_OU}
Given $0<\mu<1$, 
for a discrete-time Ornstein-Uhlenbeck (DOU) process $\{X^{\textrm{DOU}}_k, k=0, 1, \dots\}$ satisfying  
\begin{equation}
X^{\textrm{DOU}}_k=Y^{\textrm{DOU}}_k-\mu\sum_{i=0}^{k-1}{X^{\textrm{DOU}}_i}, \quad k=0, 1, \dots 
\label{eq:discrete-OU}
\end{equation} 
where $\{Y^{\textrm{DOU}}_k \}$ is a Gaussian random walk with drift $\theta$ and variance $\sigma^2$, 
the stationary density of the DOU process, $\pi$, is a normal density with mean $\theta/\mu$ and variance
$\frac{\sigma^2}{2\mu-\mu^2}$. 
\end{proposition}

\noindent \emph{Proof for Proposition~\ref{pro:stationary_density_OU}.}
Note that the DOU process $\{X^{\textrm{DOU}}_k, k=0, 1, \dots\}$ 
satisfying~(\ref{eq:discrete-OU}) is a Markov process since
\[
X^{\textrm{DOU}}_{k+1}-X^{\textrm{DOU}}_{k} = (Y^{\textrm{DOU}}_{k+1}-Y^{\textrm{DOU}}_{k})-\mu X^{\textrm{DOU}}_k. 
\]
The transition probability from state $y$ to state $x$ is
\[
\mathbb{P}(X^{\textrm{DOU}}_{k+1}=x | X^{\textrm{DOU}}_k =y) = \phi_{\theta, \sigma^2}(x-(1-\mu)y),
\]
where $\phi_{\theta, \sigma^2}(s)$ denotes the probability density function 
associated with a normal random variable with mean $\theta$ and variance $\sigma^2$.

To prove this proposition, we just need to show that 
\begin{equation}
\pi(x)=\int_{-\infty}^{\infty} \mathbb{P}(x|y)\pi(y)dy
\label{eq:discrete-OU-stationary-condition}
\end{equation}
for any given $x$, 
where
\[
\pi(x) = \frac{\sqrt{(2\mu-\mu^2)}}{\sqrt{2\pi}\sigma}\cdot
 	   \exp\left(-\frac{(2\mu-\mu^2)(x-\theta/\mu)^2}{2\sigma^2}\right). 
\]

We have 
\small
\begin{eqnarray*}
\mathbb{P} (x|y)\pi(y) &=& \frac{\sqrt{2\mu-\mu^2}}{\sqrt{2\pi}\sigma}  \frac{1}{\sqrt{2\pi}\sigma}
\exp\left(-\frac{(2\mu-\mu^2)(y-\theta/\mu)^2}{2\sigma^2}\right)
\exp\left(-\frac{(x - (1-\mu)y-\theta)^2}{2\sigma^2}\right) \\
&=& 	    \frac{\sqrt{(2\mu-\mu^2)}}{\sqrt{2\pi}\sigma} \frac{1}{\sqrt{2\pi}\sigma}
	    \exp\left(-\frac{y^2 -2[(1-\mu)x+\theta]y}{2\sigma^2} \right)
	    \exp\left(-\frac{(2\mu-\mu^2)\theta^2/\mu^2+(x-\theta)^2}{2\sigma^2} \right) \\
&=&      \frac{\sqrt{(2\mu-\mu^2)}}{\sqrt{2\pi}\sigma} \frac{1}{\sqrt{2\pi}\sigma}
	    \exp\left(-\frac{[y - ((1-\mu)x+\theta)]^2}{2\sigma^2} \right) \\
&& {}	 \cdot   \exp\left(-\frac{(2\mu-\mu^2)\theta^2/\mu^2+(x-\theta)^2-[(1-\mu)x+\theta]^2}{2\sigma^2} \right). 
\end{eqnarray*}
\normalsize 
Among which,  
\begin{eqnarray*}
V(x)      & = & \exp\left(-\frac{(2\mu-\mu^2)\theta^2/\mu^2+(x-\theta)^2-[(1-\mu)x+\theta]^2}{2\sigma^2} \right) \\
	    &=& \exp\left(-\frac{(2\mu-\mu^2)\theta^2/\mu^2+(2\mu-\mu^2)x^2-2(2-\mu)\theta x}{2\sigma^2} \right) \\
	    &=& \exp\left(-\frac{(2\mu-\mu^2)(x-\theta/\mu)^2}{2\sigma^2} \right). 
\end{eqnarray*}
Then, we have 
\begin{eqnarray*}
\int_{-\infty}^{\infty} \mathbb{P}(x|y)\pi(y)dy &=&  \frac{\sqrt{(2\mu-\mu^2)}}{\sqrt{2\pi}\sigma}
\exp\left(-\frac{(2\mu-\mu^2)(x-\theta/\mu)^2}{2\sigma^2} \right) \\
&& {} \cdot \int_{-\infty}^{\infty}  \frac{1}{\sqrt{2\pi}\sigma} \exp(-\frac{[y - ((1-\mu)x+\theta)]^2}{2\sigma^2}) dy \\
  & = & \frac{\sqrt{(2\mu-\mu^2)}}{\sqrt{2\pi}\sigma}\cdot
 	   \exp\left(-\frac{(2\mu-\mu^2)(x-\theta/\mu)^2}{2\sigma^2}\right), 
\end{eqnarray*}
which takes the exact form as the normal density with mean $\theta/\mu$ and variance $\sigma^2/(2\mu-\mu^2)$
and thus, equals to $\pi(x)$.  
This completes our proof for $\pi$ being the stationary density. 


\subsubsection{Derivation of the approximate formula}
Based on Proposition~\ref{pro:stationary_density_OU}, 
we conjecture that the stationary distribution of $X^*$ can be approximated by the following form: 
\begin{equation}
\tilde{\pi}(x)=\left\{
\begin{array}{l l}
 \pi_1(x) =  \alpha_1 \exp(-\gamma x),  & x\geq 0; \\
 \pi_2(x) = \alpha_2 \exp(-\frac{(2\mu-\mu^2)(x+\beta)^2}{2\sigma^2}), & x<0. 
\end{array}\right.
\label{eq:diffusion-density}
\end{equation}
For the ease of exposition, we use $\theta$ and $\sigma^2$
instead of $\theta_N$ and $\sigma^2_N$ to denote the mean and variance of 
the discrete-time diffusion process $X^*$. 
Moreover, in (\ref{eq:diffusion-density}), 
$\alpha_1$ and $\alpha_2$ are two normalizing constants, 
$\gamma$ is the unknown parameter for the exponential density part, 
and we define 
\[
\beta = -\mu/\theta.
\] 

The stationary density should satisfy 
\[
\tilde{\pi}(x)=\int_{-\infty}^{\infty} p(y,x) \tilde{\pi}(y) dy, 
\]
one special form of the BAR (\ref{eq:BAR}), or equivalently, 
\begin{equation}
\tilde{\pi} (x)=\int_{0}^{\infty} p(y,x)\pi_1(y) dy + \int_{-\infty}^{0} p(y,x)\pi_2(y) dy,   
\label{eq:stationary-denstity-original}
\end{equation}
where $p(y,x)$ is the transitional density of $X^*$ (from state $y$ to state $x$)
defined in (\ref{eq:transition-density-for-BAR}).

We rewrite Equation(\ref{eq:stationary-denstity-original}) as follows. 
First, for $y\geq 0$, we have
\begin{eqnarray*}
p(y,x) \pi_1(y) & = & \frac{\alpha_1}{\sqrt{2\pi}\sigma} \cdot \exp \left(-\frac{(x -y+\mu\beta)^2}{2\sigma^2} \right)
			   \cdot \exp(-\gamma y) \\
		& = & \frac{\alpha_1}{\sqrt{2\pi}\sigma}
		\cdot \exp{\left( -\frac{y^2-2(x+\mu\beta-\sigma^2\gamma)y+(x+\mu\beta)^2}{2\sigma^2} \right) } \\
		& = & \frac{\alpha_1}{\sqrt{2\pi}\sigma}
		\cdot \exp\left(-\frac{[y-(x+\mu\beta-\sigma^2\gamma)]^2}{2\sigma^2} \right)
		\cdot \exp\left(-\frac{(x+\mu\beta)^2-(x+\mu\beta-\sigma^2\gamma)^2}{2\sigma^2} \right) \\
		& = & \frac{\alpha_1}{\sqrt{2\pi}\sigma}
		\cdot \exp\left(-\frac{[y-(x+\mu\beta-\sigma^2\gamma)]^2}{2\sigma^2} \right)
		\cdot \exp\left(-\frac{\sigma^2\gamma[2x+(2\mu\beta-\sigma^2\gamma)]}{2\sigma^2} \right) \\
		& = & \frac{\alpha_1}{\sqrt{2\pi}\sigma} \cdot \exp \left(-\frac{\gamma(2\mu\beta-\sigma^2\gamma)}{2} \right)
		\cdot \exp \left(-\frac{[y-(x+\mu\beta-\sigma^2\gamma)]^2}{2\sigma^2} \right)
		\cdot \exp \left(-\gamma x \right).
\end{eqnarray*}

Therefore,
\begin{eqnarray*}
\int_{0}^{\infty} p(y,x)\pi_1(y) dy 
	&=& \alpha_1 \exp \left( -\frac{\gamma(2\mu\beta-\sigma^2\gamma)}{2} \right) \exp(-\gamma x)
	\cdot \int_{0}^{\infty} \frac{1}{\sqrt{2\pi}\sigma}\exp \left(-\frac{[y-(x+\mu\beta-\sigma^2\gamma)]^2}{2\sigma^2} \right)dy \\
	&=& \alpha_1 \exp \left(-\frac{\gamma(2\mu\beta-\sigma^2\gamma)}{2} \right) \exp(-\gamma x)
		\cdot \left[1-\Phi_{-\mu\beta, \sigma^2}\left(-x-(2\mu\beta-\sigma^2\gamma) \right) \right].
\end{eqnarray*}

Second, for $y< 0$, we have
\begin{eqnarray*}
p(y, x) \pi_2(y) & = & \frac{\alpha_2}{\sqrt{2\pi}\sigma} \cdot \exp \left(-\frac{(x-(1-\mu)y+\mu\beta)^2}{2\sigma^2} \right)
			   \cdot \exp \left(-\frac{(2\mu-\mu^2)(y+\beta)^2}{2\sigma^2} \right) \\
		& = & \frac{\alpha_2}{\sqrt{2\pi}\sigma}
		\cdot \exp \left(-\frac{(y-((1-\mu)x-\mu\beta))^2}{2\sigma^2} \right)
		\cdot \exp \left(-\frac{(2\mu-\mu^2)(x+\beta)^2}{2\sigma^2} \right).
\end{eqnarray*}

Therefore,
\begin{eqnarray*}
\int_{-\infty}^{0} p(y,x) \pi_2(y) dy & = & \int_{-\infty}^{0} \frac{\alpha_2}{\sqrt{2\pi}\sigma}
		\cdot \exp \left(-\frac{(y-((1-\mu)x-\mu\beta))^2}{2\sigma^2} \right)
		\cdot \exp \left(-\frac{(2\mu-\mu^2)(x+\beta)^2}{2\sigma^2} \right) dy \\
	&=& \alpha_2 \exp \left(-\frac{(2\mu-\mu^2)(x+\beta)^2}{2\sigma^2} \right)
	\cdot \int_{-\infty}^{0} \frac{1}{\sqrt{2\pi}\sigma}\exp \left(-\frac{(y-(1-\mu)x +\mu\beta)^2}{2\sigma^2} \right)dy \\
	&=& \alpha_2 \exp \left(-\frac{(2\mu-\mu^2)(x+\beta)^2}{2\sigma^2} \right) 
	\cdot \Phi_{-\mu\beta, \sigma^2}(-x+\mu x).
\end{eqnarray*}


If (\ref{eq:stationary-denstity-original}) holds, when $x\geq 0$, 
we should have 
\begin{eqnarray*}
\alpha_1 \exp(-\gamma x) &=& \alpha_1 \exp \left(-\frac{\gamma(2\mu\beta-\sigma^2\gamma)}{2} \right) \exp(-\gamma x) 
\cdot \left[1-\Phi_{-\mu\beta, \sigma^2}\left(-x-(2\mu\beta-\sigma^2\gamma) \right) \right] \\
&& {} +\alpha_2 \exp \left(-\frac{(2\mu-\mu^2)(x+\beta)^2}{2\sigma^2} \right) 
\cdot \Phi_{-\mu\beta, \sigma^2}(-x+\mu x), 
\end{eqnarray*}
which is equivalent to
\begin{eqnarray}
&& \alpha_1 \exp(-\gamma x) \cdot
\left[1-\exp \left(-\frac{\gamma(2\mu\beta-\sigma^2\gamma)}{2} \right) 
\left(1- \Phi_{-\mu\beta, \sigma^2}(-x-(2\mu\beta-\sigma^2\gamma)) \right) \right]  \nonumber \\ 
&=&  \alpha_2 \exp \left(-\frac{(2\mu-\mu^2)(x+\beta)^2}{2\sigma^2} \right) 
\cdot \Phi_{-\mu\beta, \sigma^2}(-x+\mu x).
\label{eq:to-verify-1}
\end{eqnarray}
\normalsize

Similarly, if (\ref{eq:stationary-denstity-original}) holds, 
when $x<0$, we should have 
\begin{eqnarray*}
\alpha_2 \exp \left(-\frac{(2\mu-\mu^2)(x+\beta)^2}{2\sigma^2} \right)
&=&  \alpha_1 \exp \left(-\frac{\gamma(2\mu\beta-\sigma^2\gamma)}{2} \right) \exp(-\gamma x)
		\cdot \left[1-\Phi_{-\mu\beta, \sigma^2}\left(-x-(2\mu\beta-\sigma^2\gamma) \right) \right]  \\
&& {} +\alpha_2 \exp \left(-\frac{(2\mu-\mu^2)(x+\beta)^2}{2\sigma^2} \right) 
\cdot \Phi_{-\mu\beta, \sigma^2}(-x+\mu x),
\end{eqnarray*}
which is equivalent to
\begin{eqnarray}
&& \alpha_1 \exp \left(-\frac{\gamma(2\mu\beta-\sigma^2\gamma)}{2} \right) \exp(-\gamma x)
		\cdot \left[1-\Phi_{-\mu\beta, \sigma^2}(-x-(2\mu\beta-\sigma^2\gamma)) \right] \nonumber \\
&=& \alpha_2 \exp \left(-\frac{(2\mu-\mu^2)(x+\beta)^2}{2\sigma^2} \right) 
\cdot \left(1-\Phi_{-\mu\beta, \sigma^2}(-x+\mu x) \right).
\label{eq:to-verify-2}
\end{eqnarray}

When $x=0$, Equations (\ref{eq:to-verify-1}) and (\ref{eq:to-verify-2}) become  
\begin{equation}
\alpha_1 \cdot
\left[1-\exp \left(-\frac{\gamma(2\mu\beta-\sigma^2\gamma)}{2} \right)
\left(1- \Phi_{-\mu\beta, \sigma^2}(-(2\mu\beta-\sigma^2\gamma)) \right) \right] 
= \alpha_2 \exp \left(-\frac{(2\mu-\mu^2)\beta^2}{2\sigma^2} \right) \cdot \Phi_{-\mu\beta, \sigma^2}(0),
\label{eq:verify1-x=0}
\end{equation}
and
\begin{equation}
\alpha_1 \exp \left(-\frac{\gamma(2\mu\beta-\sigma^2\gamma)}{2} \right)
	      \cdot \left[1-\Phi_{-\mu\beta, \sigma^2}(-(2\mu\beta-\sigma^2\gamma)) \right]
= \alpha_2 \exp \left(-\frac{(2\mu-\mu^2)\beta^2}{2\sigma^2} \right) \cdot \left(1-\Phi_{-\mu\beta, \sigma^2}(0) \right), 
\label{eq:verify2-x=0}
\end{equation}
respectively. 

Recall that $\pi(x)$ is continuous at $x=0$. Thus, the two normalizing constants satisfy: 
\begin{equation}
\pi_1(0)=\alpha_1 = \alpha_2 \exp \left(-\frac{(2\mu-\mu^2)\beta^2}{2\sigma^2} \right)=\pi_2(0).
\label{eq:a1-a2-relation1}
\end{equation}

Comparing (\ref{eq:a1-a2-relation1}) with
(\ref{eq:verify1-x=0}) and (\ref{eq:verify2-x=0}), we find that when 
\[
\sigma^2\gamma = 2\mu\beta = - 2 \theta,
\]
or equivalently, 
\begin{equation}
\gamma = -\frac{2\theta}{\sigma^2},
\label{gamma} 
\end{equation}
both (\ref{eq:verify1-x=0}) and (\ref{eq:verify2-x=0}) can be satisfied.
Therefore, we choose $\gamma$ in (\ref{gamma}), which eventually gives us (\ref{eq:diffusion-density-revised}). 

Unfortunately, using this $\gamma$, 
we are unable to show (\ref{eq:to-verify-1}) and (\ref{eq:to-verify-2}) hold for a general $x$.


\section{Numerical results on diffusion approximations}
\label{sec:numerical-results-for-diffusion}

\subsection{Approximation for the midnight count distribution}

\begin{figure}[btp]
        \centering
        \subfloat[$N=500, \Lambda=90.95, \rho=0.96$]{
       \includegraphics[width=0.31\linewidth, height = 0.15\textheight]
        {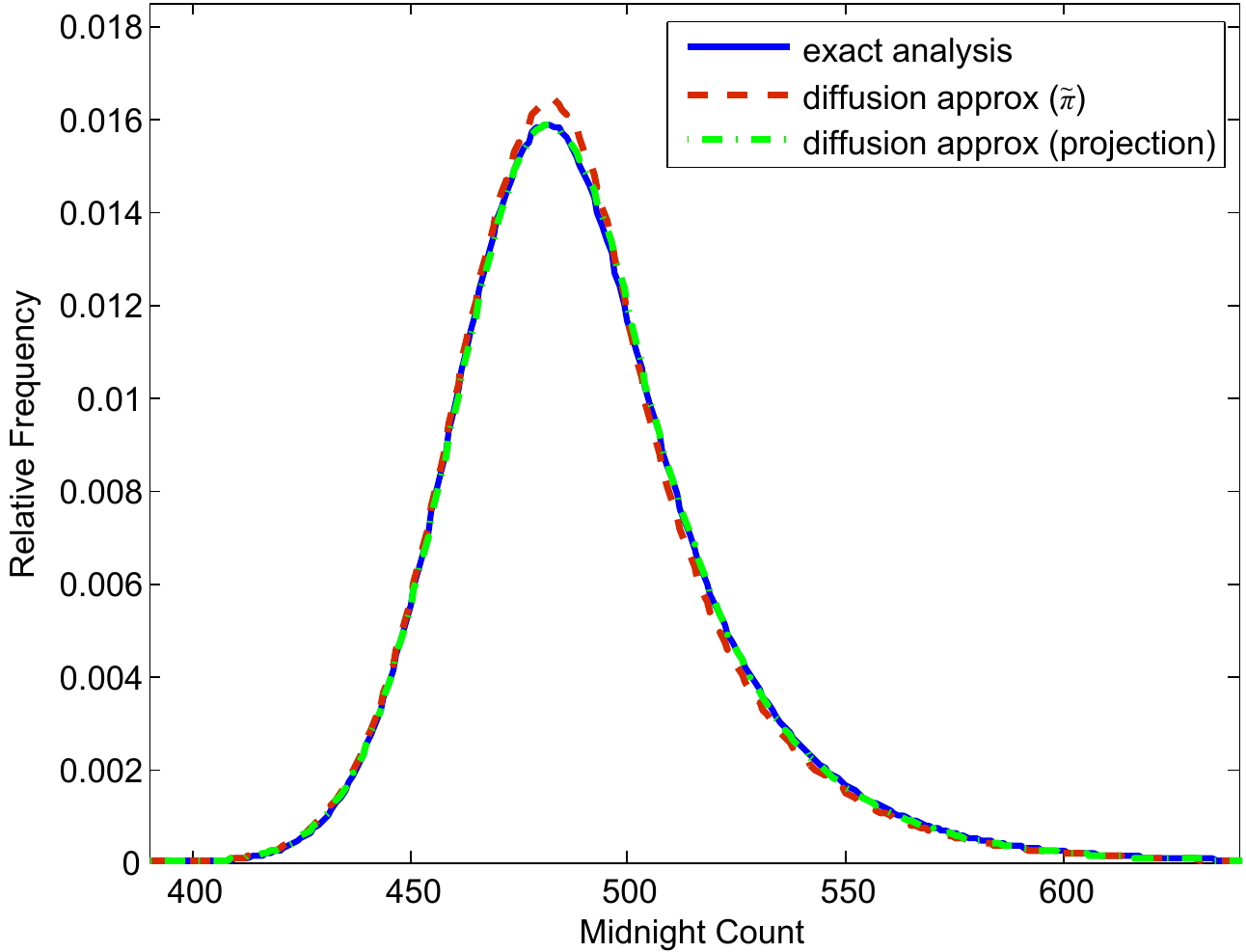}
         \label{fig:mignight-compare-N-500} }
        \ 
       \subfloat[$N=66, \Lambda=11.37, \rho=0.91$]{
       \includegraphics[width=0.31\linewidth, height = 0.15\textheight]
       {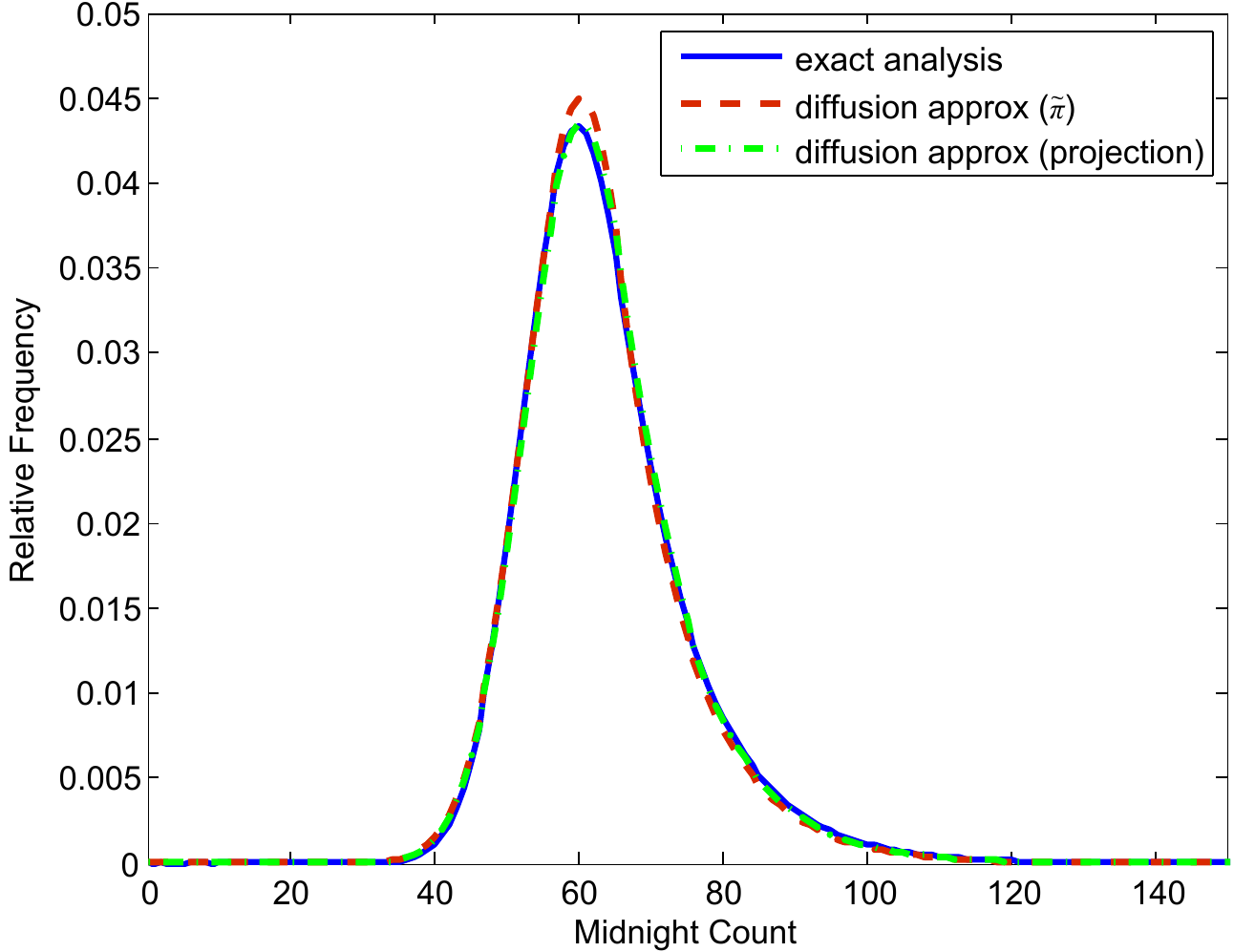}
         \label{fig:mignight-compare-N-66} }
         \ 
       \subfloat[$N=18, \Lambda=3.03, \rho=0.89$]{
       \includegraphics[width=0.31\linewidth, height = 0.15\textheight]
         {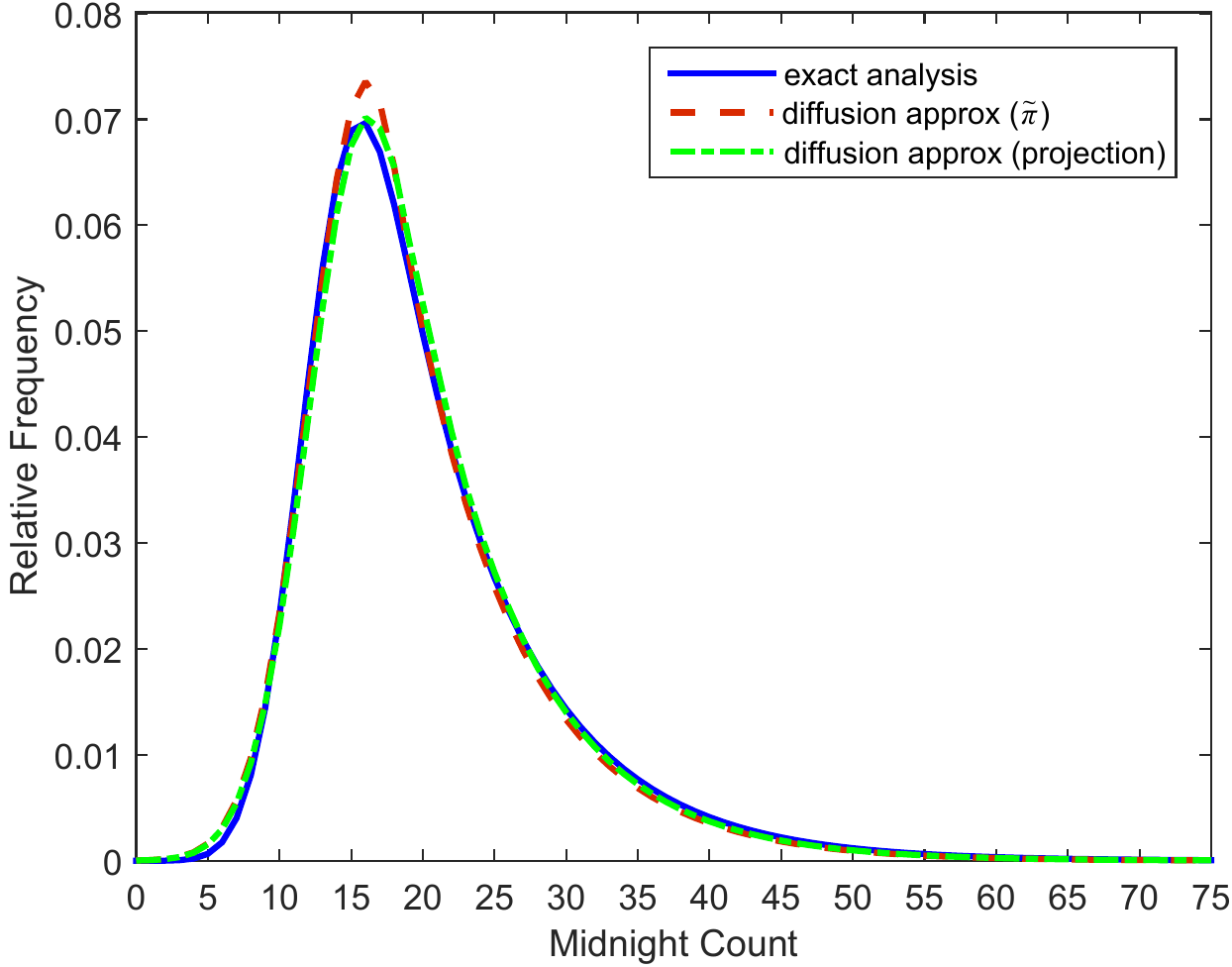}
         \label{fig:mignight-compare-N-18} }
 \caption[Stationary distribution of the midnight customer count from exact Markov chain analysis
 and diffusion approximations.]
 {\small{Stationary distribution of the midnight customer count from exact Markov chain analysis
 and diffusion approximations.
 Here, the mean LOS is 5.3 days,
 and we do not specify the discharge distribution 
 because it does not affect the midnight customer count distribution.}}
 \label{fig:numerical-comp-diffusion-density-midnight}
 \end{figure}

Figure~\ref{fig:numerical-comp-diffusion-density-midnight}
compares the stationary distributions of the midnight customer count
solved (i) from the exact Markov chain analysis,
(ii) from using the approximate formula $\tilde{\pi}$ in (\ref{eq:diffusion-density-revised}),
and (iii) from using the projection algorithm specified in Section~\ref{sec:projection-algorithm}.
The parameter settings for these numerical experiments are the same as those in Section~5 of the main paper.
We test a large system ($N=500$) and two small systems ($N=66$ and $18$),
with the utilization $\rho$ being 96\%, 91\% and 89\%, respectively.


\subsection{Time-dependent performance}

\begin{figure}[btp]
\centering
        \subfloat[$\E_\infty (Q(t)),\ N=500$]
        {
                \centering
                \includegraphics[width=0.31\linewidth, height = 0.15\textheight]
                {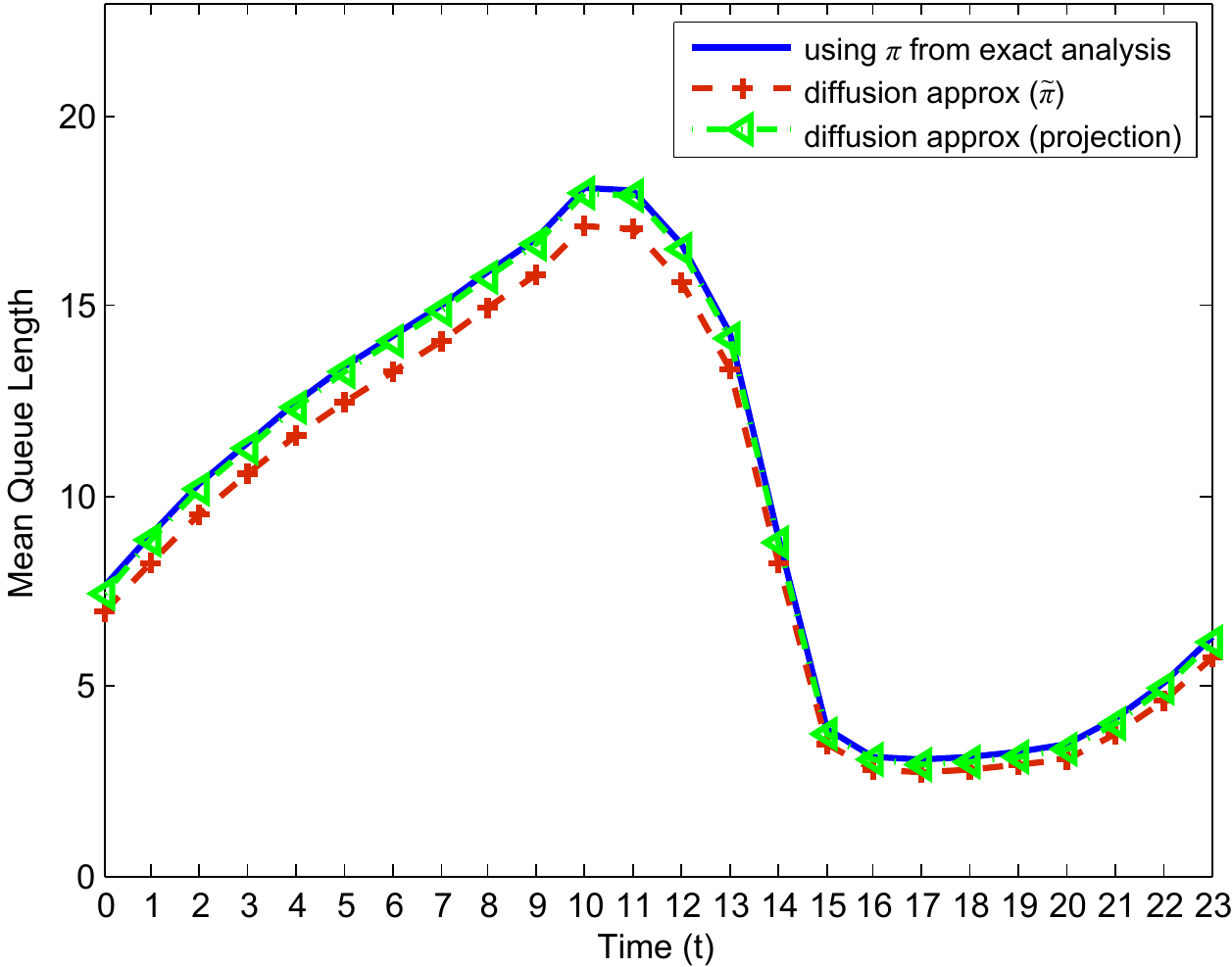}
                }
        \
        \subfloat[$\E_\infty(W(t)),\ N=500$]
        {
                \centering
                \includegraphics[width=0.31\linewidth, height = 0.15\textheight]
                {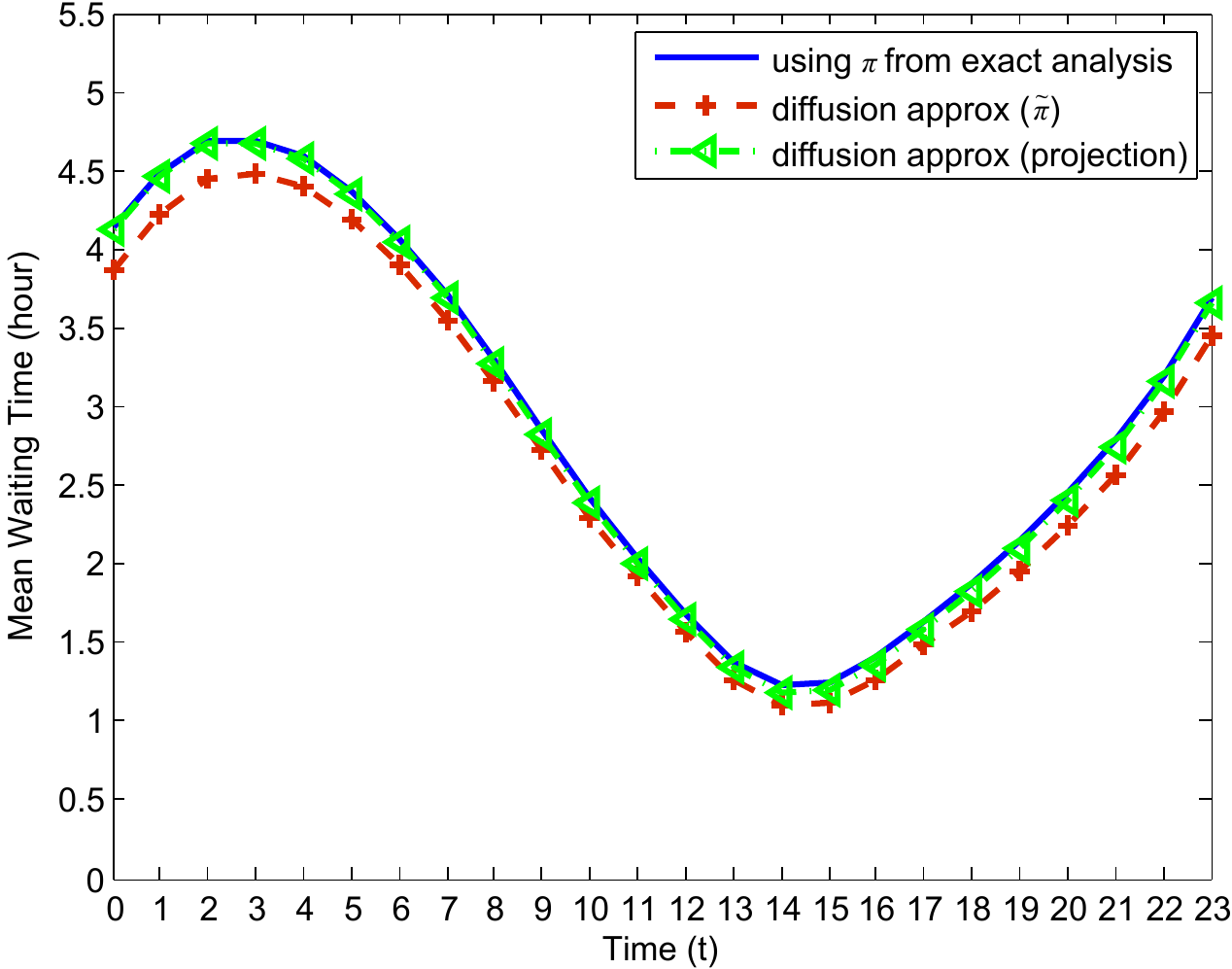}
        }
        \
        \subfloat[$P_\infty(W(t)>6/24),\ N=500$]
        {
                \centering
                \includegraphics[width=0.31\linewidth, height = 0.15\textheight]
                {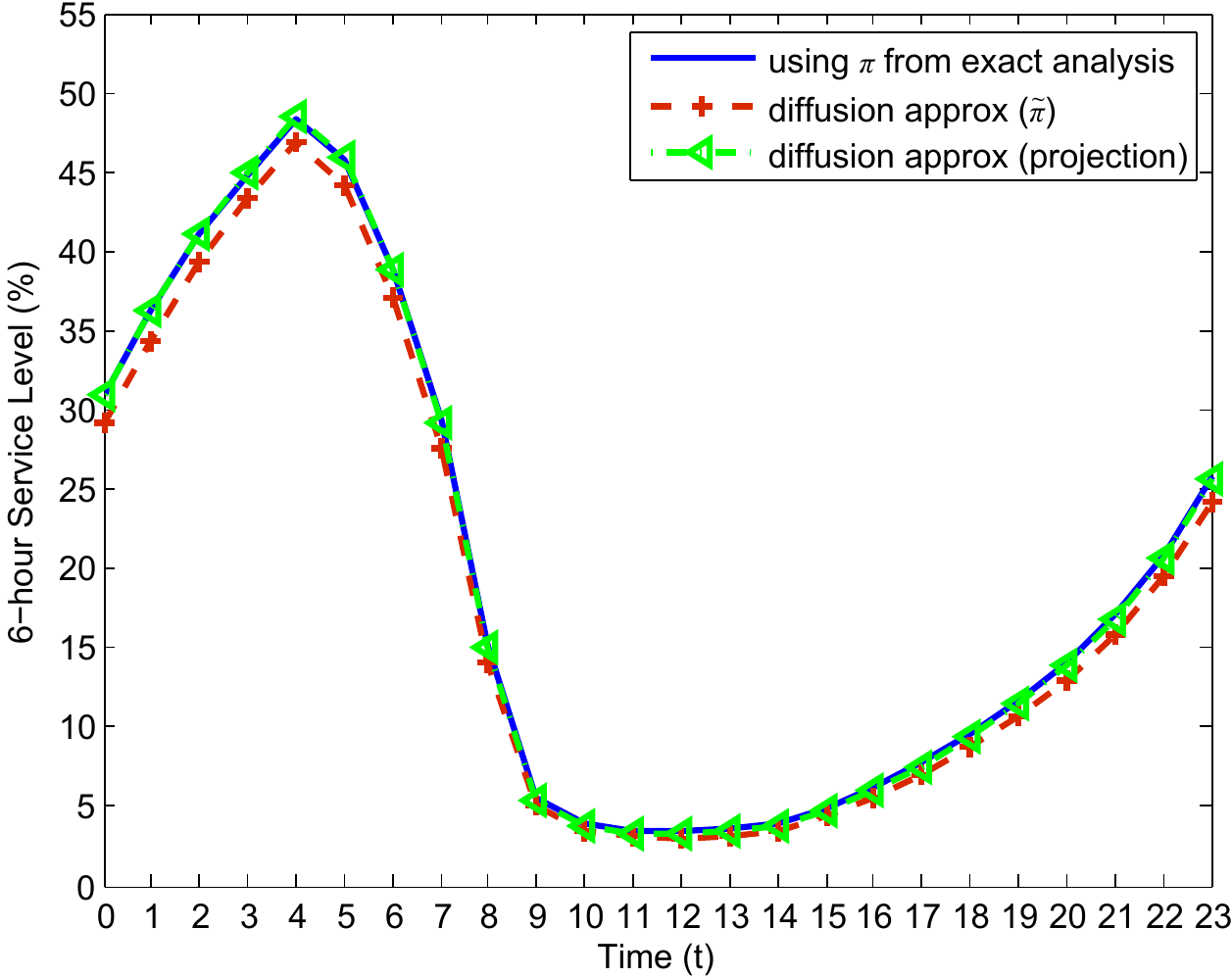}
        }
\caption[Time-dependent performance curves from exact analysis
 and diffusion approximations.]
 {\small{Time-dependent performance curves from exact analysis
 and diffusion approximations.
 Here, $\Lambda=90.95$ for $N=500$. 
 We fix the mean LOS as 5.3 days and use the baseline discharge distribution.
 The three performance curves in each subfigure are from normal approximations
 using (i) $\pi$ solved from exact Markov chain analysis,
 (ii) $\tilde{\pi}$ in (\ref{eq:diffusion-density-revised}), 
and (iii) $\pi^*$ solved from the projection algorithm, respectively.}}
\label{fig:hourly-waiting-exact-approx-500-gen-dis}
\end{figure}

\begin{figure}[btp]
\centering
	\subfloat[$\E_\infty (Q(t)),\ N=66$]
        {
                \centering
                \includegraphics[width=0.31\linewidth, height = 0.15\textheight]
                {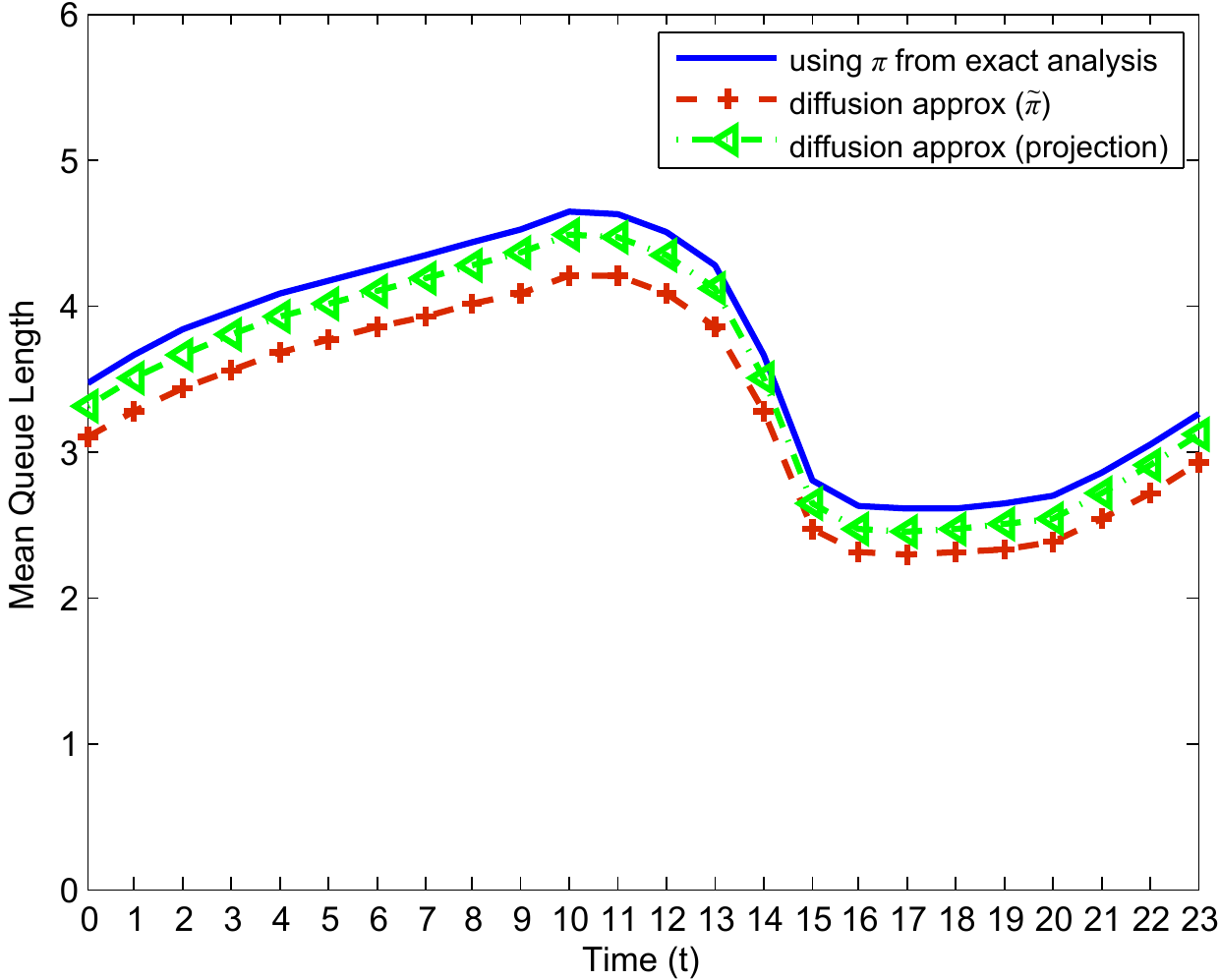}
                }
        \
        \subfloat[$\E_\infty(W(t)),\ N=66$]
        {
                \centering
                \includegraphics[width=0.31\linewidth, height = 0.15\textheight]
                {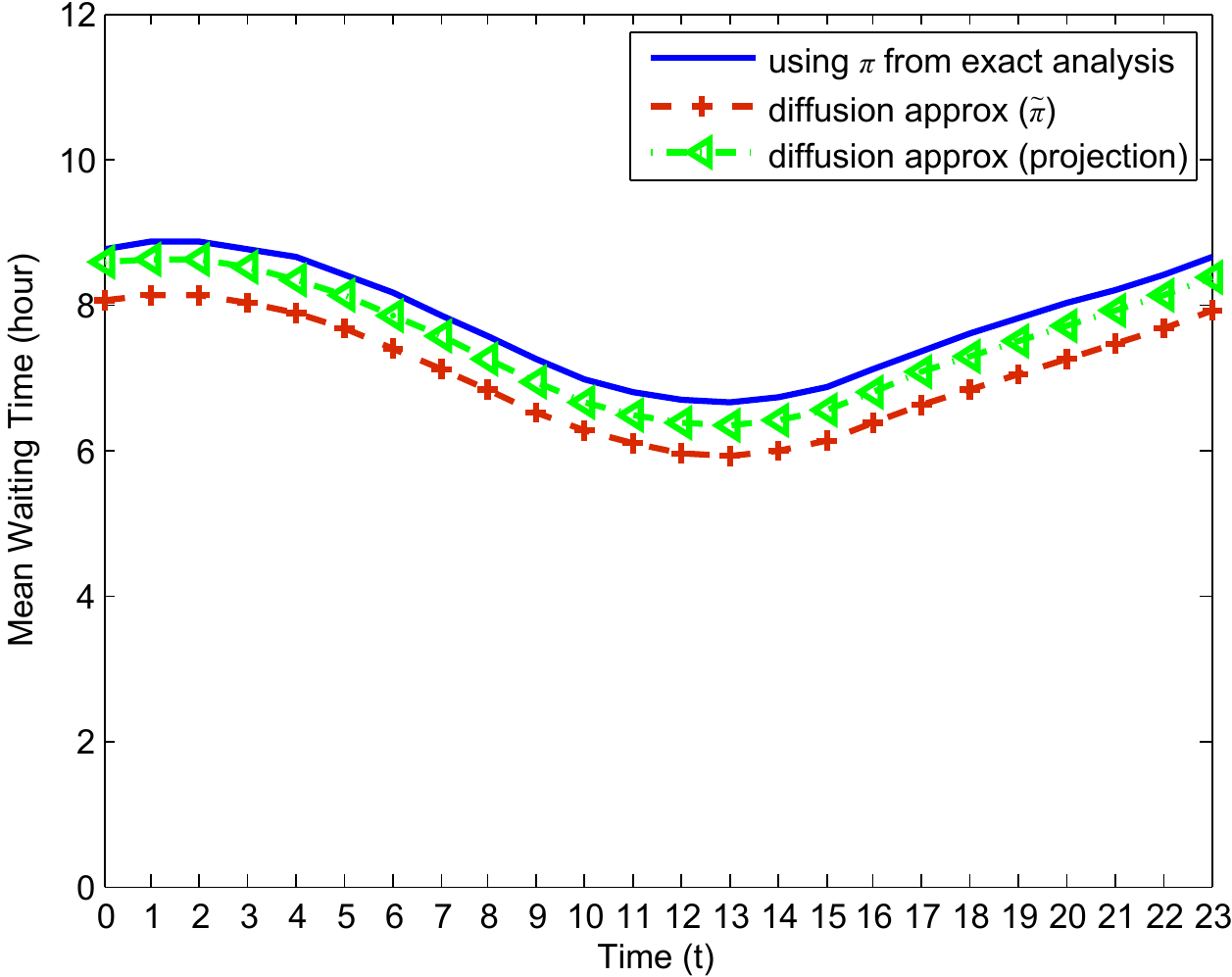}
        }
        \
        \subfloat[$P_\infty(W(t)>6/24),\ N=66$]
        {
                \centering
                \includegraphics[width=0.31\linewidth, height = 0.15\textheight]
                {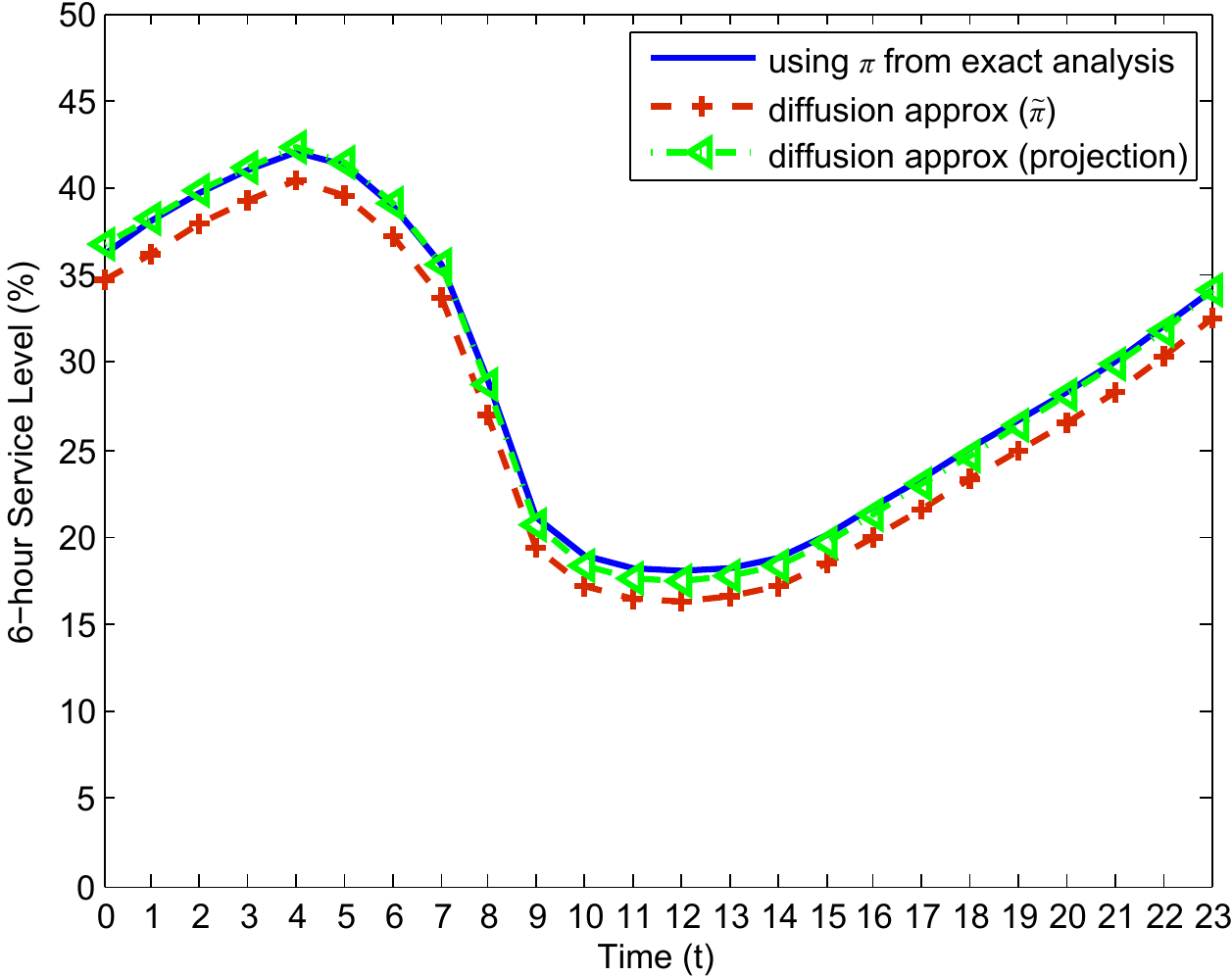}
        }
\caption[Time-dependent performance curves from exact analysis
 and diffusion approximations.]
 {\small{Time-dependent performance curves from exact analysis
 and diffusion approximations.
 Here, $\Lambda=11.37$ for $N=66$.
 We fix the mean LOS as 5.3 days and use the baseline discharge distribution. 
 The three performance curves in each subfigure are from normal approximations
 using (i) $\pi$ solved from exact Markov chain analysis,
 (ii) $\tilde{\pi}$ in (\ref{eq:diffusion-density-revised}), 
and (iii) $\pi^*$ solved from the projection algorithm, respectively.}}
\label{fig:hourly-waiting-exact-approx-66-gen-dis}
\end{figure}

\begin{figure}[btp]
\centering
	\subfloat[$\E_\infty (Q(t)),\ N=18$]
        {
                \centering
                \includegraphics[width=0.31\linewidth, height = 0.15\textheight]
                {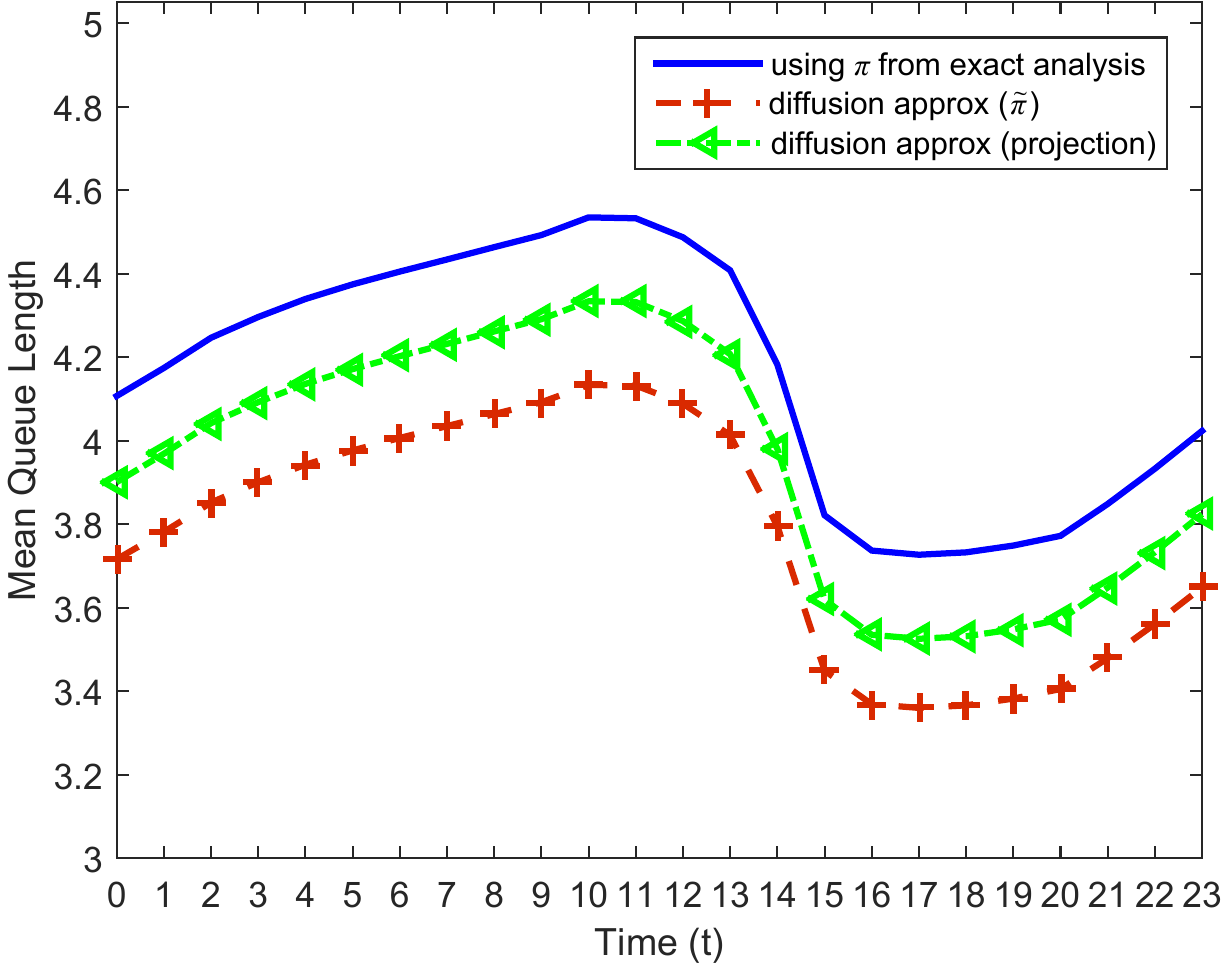}
                }
        \
        \subfloat[$\E_\infty(W(t)),\ N=18$]
        {
                \centering
                \includegraphics[width=0.31\linewidth, height = 0.15\textheight]
                {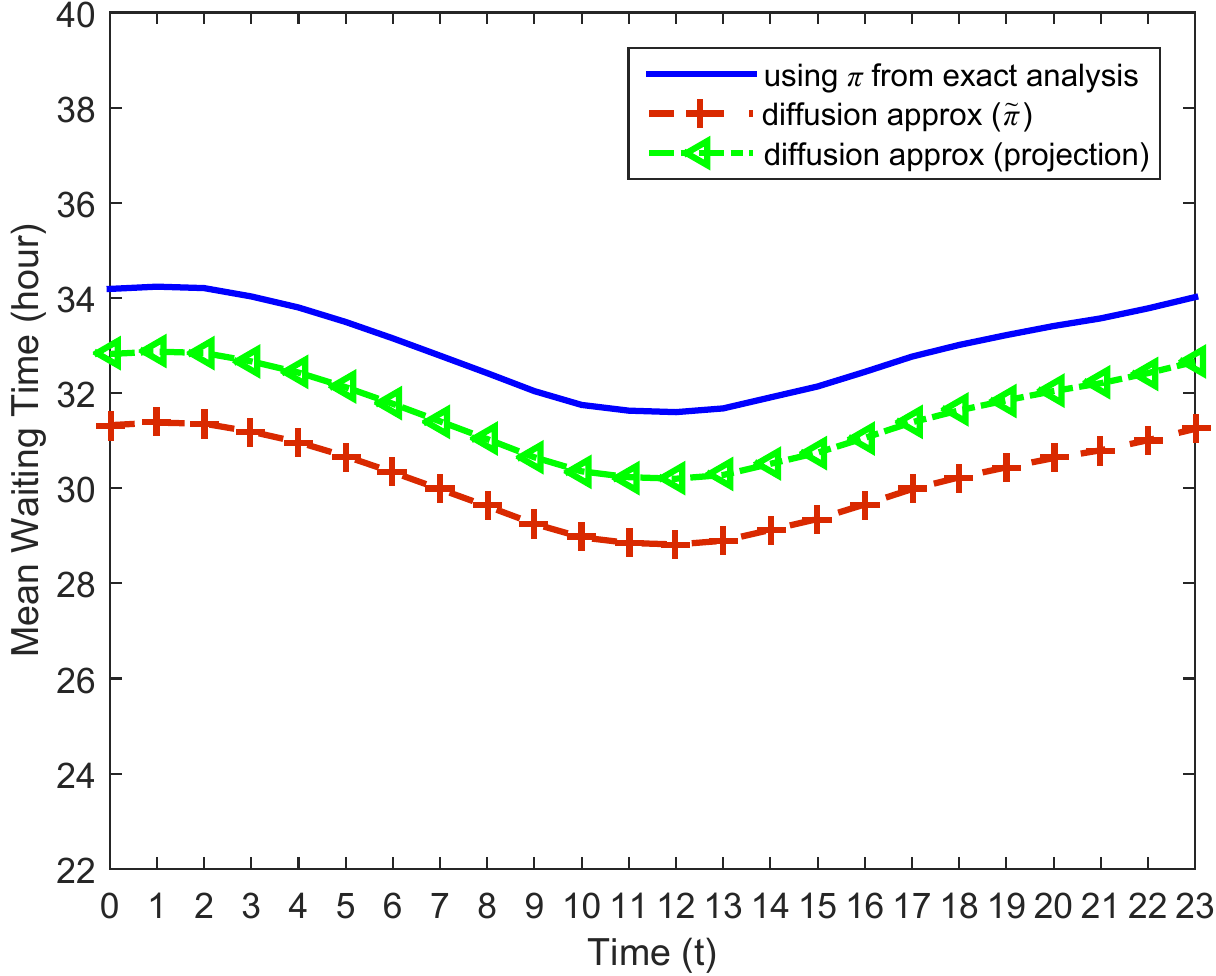}
        }
        \
        \subfloat[$P_\infty(W(t)>6/24),\ N=18$]
        {
                \centering
                \includegraphics[width=0.31\linewidth, height = 0.15\textheight]
                {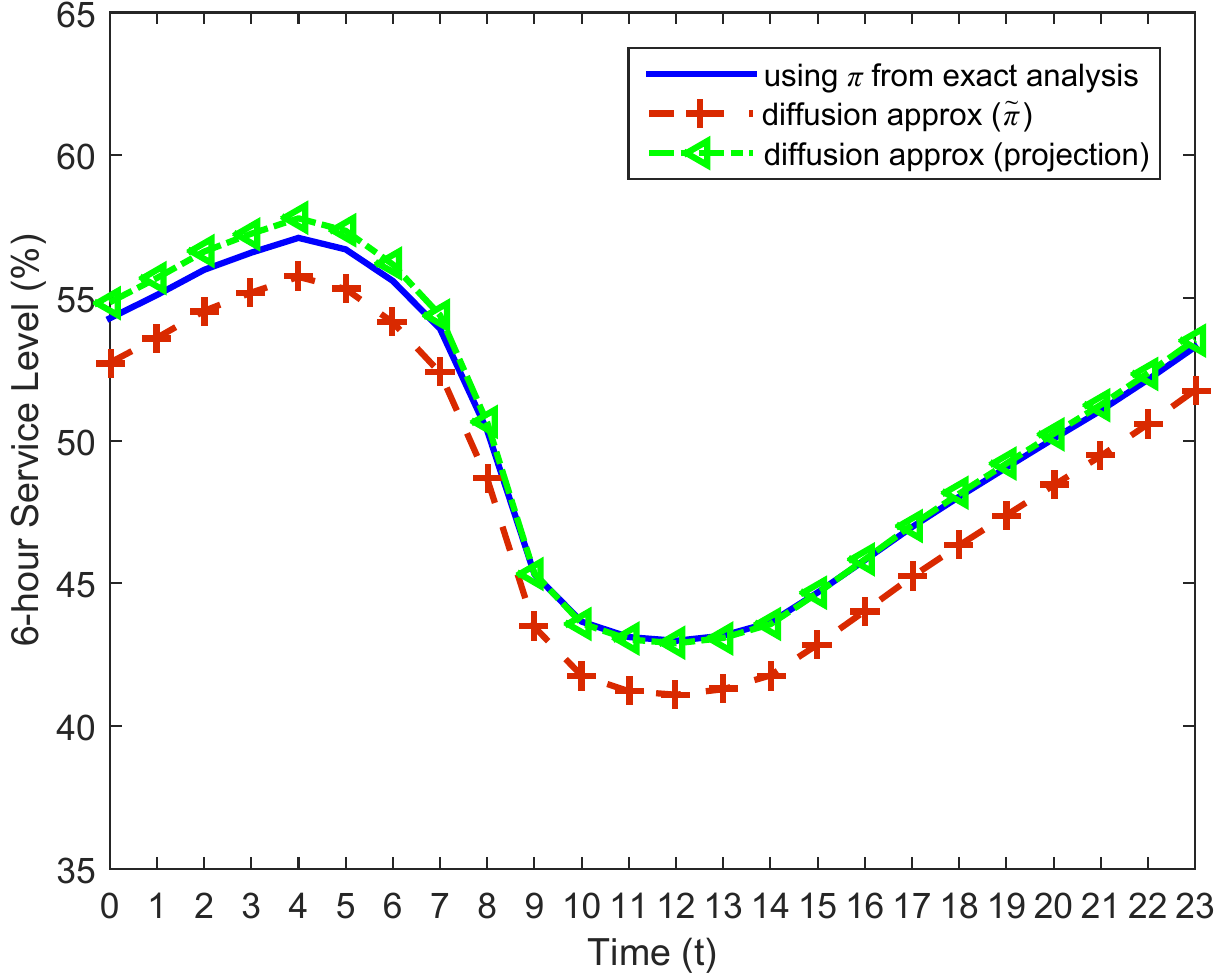}
        }
\caption[Time-dependent performance curves from exact analysis
 and diffusion approximations.]
 {\small{Time-dependent performance curves from exact analysis
 and diffusion approximations.
 Here, $\Lambda=3.03$ for $N=18$.
 We fix the mean LOS as 5.3 days and use the baseline discharge distribution.
 The three performance curves in each subfigure are from normal approximations
 using (i) $\pi$ solved from exact Markov chain analysis,
 (ii) $\tilde{\pi}$ in (\ref{eq:diffusion-density-revised}), 
 and (iii) $\pi^*$ solved from the projection algorithm, respectively.}}
\label{fig:hourly-waiting-exact-approx-18-gen-dis}
\end{figure}

Figures~\ref{fig:hourly-waiting-exact-approx-500-gen-dis}, 
~\ref{fig:hourly-waiting-exact-approx-66-gen-dis}
and~\ref{fig:hourly-waiting-exact-approx-18-gen-dis}
show the time-dependent performance for systems with $N=500,\ 66$, and $18$, respectively. 
The three curves in each subfigure are obtained from normal approximations
using (i) $\pi$ solved from exact Markov chain analysis,
(ii) $\tilde{\pi}$ in (\ref{eq:diffusion-density-revised}), 
and (iii) $\pi^*$ solved from the projection algorithm specified in Section~\ref{sec:projection-algorithm}. 

\bibliography{dai10012012,nuh}

\end{document}